\documentclass[12pt]{amsart}
\usepackage{anysize}
\usepackage{amsmath}
\usepackage{amsthm}
\usepackage{amsmath,amscd}
\usepackage{amssymb}
\usepackage{stmaryrd}
\usepackage{wasysym}
\usepackage{mathrsfs}
\usepackage{hyperref}
\usepackage{dsfont}
\usepackage{soul}
\hypersetup{colorlinks=true,linkcolor=blue,citecolor=black}
\usepackage{color}
\usepackage{float}
\usepackage{bm}
\usepackage{algorithm}
\usepackage[noend]{algpseudocode}
\usepackage{mathtools}
\usepackage{comment}
\usepackage{thmtools}
\usepackage{enumerate}
\usepackage{graphicx}
\graphicspath{ {./Figures/} }
\usepackage{setspace}
\newcommand\mtop{1in}
\newcommand\mbottom{1in}
\newcommand\mleft{1in}
\newcommand\mright{1in}
\usepackage[top = \mtop, bottom = \mbottom, left = \mleft, right=\mright]{geometry}
\usepackage{subcaption}

\newtheorem{Theorem}{Theorem}[section]
\newtheorem{corollary}{Corollary}[Theorem]
\newtheorem{lemma}[Theorem]{Lemma}

\newtheorem{proposition}{Proposition}[section]

\theoremstyle{definition}
\newtheorem{definition}{Definition}

\newcommand{\R}{\mathbb{R}}
\newcommand{\Z}{\mathbb{Z}}

\def\P{{\mathbb{P}}}

\usepackage{scalerel,stackengine}
\stackMath
\newcommand\reallywidehat[1]{%
	\savestack{\tmpbox}{\stretchto{%
			\scaleto{%
				\scalerel*[\widthof{\ensuremath{#1}}]{\kern-.6pt\bigwedge\kern-.6pt}%
				{\rule[-\textheight/2]{1ex}{\textheight}}
			}{\textheight}%
		}{0.5ex}}%
	\stackon[1pt]{#1}{\tmpbox}%
}

\title{Equivalence of Polychromatic Arm Probabilities on the Square Lattice}

\author{Lily Reeves}
\address{Center for Applied Mathematics, Cornell University. Frank H.T. Rhodes Hall, Ithaca 14850.}
\email{zw477@cornell.edu}

\author{Philippe Sosoe}
\thanks{P.S.'s research is partially supported by NSF grant DMS 1811093.}
\address{Department of Mathematics, Cornell University. Malott Hall, Ithaca, NY 14853.}
\email{psosoe@math.cornell.edu}

\begin{document}
	\maketitle
	
	\begin{abstract}
		We consider $2d$ critical Bernoulli percolation on the square lattice. We prove an approximate color-switching lemma  comparing  $k$-arm probabilities for different polychromatic color sequences. This result is well-known for site percolation on the triangular lattice in \cite{Nolin08}. To handle the complications arising from the dual lattice, we introduce a shifting transformation to convert arms between the primal and dual lattices.
	\end{abstract}

\section{Introduction}

We consider Bernoulli percolation on the square lattice $\Z^2$: each edge is colored open or closed with probability $p$ and $1-p$, respectively. We select the critical parameter $p=p_c=\frac{1}{2}$. 
In this setting, \emph{arm events}, which are defined by the simultaneous occurrence of several long-range connections across annuli, have been extensively studied in critical percolation, especially since the pioneering works by Harry Kesten in the 1980s, see \cite{Kesten80, Kesten86, Kesten87, SW01}. See \cite{Nolin08} for a survey. 

Configurations in a $k$-arm event have $k$ disjoint paths, each of a single color (open or closed-dual), connecting the two boundaries of an annulus. A $k$-arm event is said to be \emph{polychromatic} if not all paths in the configuration are of the same color. It is well-known that for critical \emph{site} percolation on the triangular lattice, probabilities of different polychromatic arm events are comparable up to constants that are independent of the size of the annulus. In this note, we prove that the same holds on the square lattice.

Some of the key results accessible for site percolation on the triangular lattice remain out of reach for the square lattice. 
While many of the techniques involving the use of correlation inequalities, e.g. generalized FKG inequality and BK inequality, can be applied to very general lattices, intrinsic differences in the duality relations satisfied by the two models have so far impeded the extension of key results for critical and near-critical percolation on the triangular lattice to bond percolation on the square lattice. For example, the arm exponents for critical site percolation on the triangular lattice are known \cite{LSW02,SW01}: $\alpha_1 = 5/48$ and $\alpha_k = (k^2-1)/12$ for $k\geq 2$. Since this result relies on approximation by Schramm-Loewner Evolution (SLE), the values of arm exponents on the square lattice remain conjectural. 

The proof of the equivalence of different polychromatic arm probabilities on the triangular lattice uses \emph{color switching}, a combinatorial trick with a venerable history in critical percolation, see e.g. \cite{ADA99}. Most notably, an \emph{exact} version of color-switching is used in Smirnov's celebrated proof of conformal invariance for critical site percolation \cite{Smirnov01}. 

We adapt to the square lattice an approximate color-switching argument in \cite{Nolin08} which shows that arm probabilities for different color sequences of the same length are asymptotically equivalent as long as they contain arms of both colors: 
\begin{proposition}\label{lemma:mainresult}
    Consider percolation on the square lattice and let $k\geq 3$, $n_0(k)< n < N$. Let $\sigma, \sigma'$ be two polychromatic color sequences (consisting of open and dual-closed  arms) of length $k$. We denote $A_{k,\sigma}(n,N)$ to be the event that there exist $k$ arms, color-coded by $\sigma$, from scale $n$ to $N$. See the notation section for precise definitions. Then, there exists a constant $C$ independent of choices of $k,\sigma, \sigma',n$, and $N$, such that
    \begin{equation*}
        \P(A_{k,\sigma}(n,N)) \leq C \P(A_{k,\sigma'}(n,N)).
    \end{equation*}
\end{proposition}

\subsection{Notations} \label{subsection:notations}

In this section, we summarize the notations we will use. Throughout the note, we consider Bernoulli percolation on the square lattice $\Z^2$ seen as a graph with the edge set $\mathcal{E}$ consisting of all pairs of nearest-neighbor vertices. 

We let $\P$ be the critical bond percolation measure 
\[\P = \prod_{e\in \mathcal{E}} \frac12 (\delta_0+\delta_1)\] 
on the state space $\Omega = \{0,1\}^{\mathcal{E}}$, with the product $\sigma$-algebra. An edge $e$ is said to be open in the configuration $\omega\in \Omega$ if $\omega(e) = 1$ and closed otherwise. 

A (lattice) \textit{path} is a sequence $(v_0, e_1, v_1, \dots , v_{N-1}, e_N , v_N)$ such that for all $k = 1,\dots,N$, $\|v_{k-1} - v_k\|_1 = 1$ and $e_k = \{v_{k-1},v_k\}$. 
Given $\omega \in \Omega$, we say that $\gamma=(e_k)_{k=1,\dots,N}$ is open in $\omega$ if $\omega(e_k) = 1$ for $k = 1,\dots ,N$.



The dual lattice is written $((\Z^2)^*, \mathcal{E}^*)$, where $\left(\Z^2\right)^* = \Z^2 + (\frac12, \frac12)$
with its nearest-neighbor edges.  Given $\omega \in \Omega$, we obtain a configuration $\omega^* \in \Omega^* = \{0,1\}^{\mathcal{E}^*}$ by the relation $\omega^*(e^*) = \omega(e)$, where $e^*$ is the dual edge that shares a midpoint with $e$. For any $V \subset \mathbb{R}^2$ we write $V^*=V+(\frac{1}{2},\frac{1}{2}).$

For $x \in \Z^2$, we define 
\[B(x,n)=\{(y, z)\in \mathcal{E} : \|y-x\|_\infty \le n, \|z-x\|_\infty \le n\}.\]
Here $x\sim y$ means $x$ and $y$ are nearest neighbors on the lattice $\Z^2$, and $\|x-y\|_\infty = \max_{i=1,2}|x_i-y_i|$.
When $x$ is the origin $(0,0)$, we sometimes abbreviate $B((0,0),n)$ by $B_n$ or $B(n)$. We denote by $\partial B(x,n)$ the set 
\[\partial B(x,n)=\{(x_1, x_2)\in \mathcal{E} : \|x_1-x\|_\infty = n, \|x_2-x\|_\infty = n\}.\]
We sometimes say that a vertex $x$ \emph{lies in} $\partial B(x,n)$ if this vertex coincides with an endpoint of an edge in $\partial B(x,n)$.
We define an annulus centered at $x \in \Z^2$ as the difference between two boxes of different sizes centered at $x$:
\[\text{For $0<n<N$, }  B(x,n,N) = B(x,N) \setminus B(x,n).\]
We often abbreviate $B((0,0),n,N)$ as $B(n,N)$ when $x=(0,0)$ is implied.




A \emph{color sequence} $\sigma$ of length $k$ is a sequence $(\sigma_1, \dots, \sigma_k)\in \{O, O^*, C, C^*\}^k$. Each $\sigma_i$ indicates a ``color''. The colors are encoded $O$ for open, $O^*$ for dual-open, $C$ for closed, and $C^*$ for dual-closed. From the percolation point of view, the most interesting color sequences consist of primal-open and/or dual-closed connections since by duality, primal-open paths and dual-closed paths are mutually disjoint. But in the proof of the main Proposition, it will be useful to consider primal-closed and dual-open arms, as well. 

We use the convention $(O^*)^*=O$, $(C^*)^*=C$, and $\sigma^* = (\sigma_1, \dots, \sigma_k)^* = (\sigma_1^*, \dots, \sigma_k^*)$. Similarly, $\overline{\sigma}$ denotes the flipped color sequence, with the conventions $\overline{O}=C$, $\overline{C}=O$, $\overline{O^*}=C^*$, $\overline{C^*}=O^*$, and $\overline{(\sigma_1, \dots, \sigma_k)}$ = ($\overline{\sigma_1}, \dots, \overline{\sigma_k}$).

\subsection*{Arm events}
A primal open (closed, resp.) arm in $B(n,N)$ connecting $\partial B_n$ and $\partial B_N$ is a path of open (closed, resp.) edges in $B(n,N)$ with one endpoint lying in $\partial B_n$ and another endpoint in $\partial B_N$. 

We say a \textit{dual} open (closed, resp.) arm in $B(n,N)^*$ \emph{connects} $\partial B_n$ and $\partial B_N$ if the (primal) path obtained by shifting by $(-1/2,-1/2)$ connects $\partial B_n$ and $\partial B_N$.

\begin{definition}	For $n\leq N$, we define the $k$-arm event with color sequence $\sigma$ to be the event that there are $k$ disjoint paths whose colors are specified by $\sigma$ in clockwise order in the annulus $B(n,N)$ connecting $\partial B_n$ and $\partial B_N$. Formally,
    \begin{equation*}
        A_{k,\sigma}(n,N) := \{\partial B_n \leftrightarrow_{k,\sigma} \partial B_N\}.
    \end{equation*}
\end{definition}

We note a technical point: for $A_{k,\sigma}(n,N)$ to be defined, $n$ needs to be large enough for all $k$ arms to be (vertex)-disjoint. We define $n_0(k)$ to be the smallest integer such that $|\partial B(n_0(k))| \geq k$. Color sequences that are equivalent up to cyclic order denote the same arm event.

\subsection{Outline of the paper}
We use the basic idea of color switching. Unlike in site percolation on the triangular lattice, when flipping the statuses of edges in a region, a primal-open arm becomes primal-closed instead dual-closed; a dual-closed arm becomes dual-open instead of primal-open. To address this problem, we introduce a shifting transformation in Section \ref{section:shift-transformation} to convert between the two lattices. To apply the transformation, the arms cannot come too close to the boundaries. Thus, in Section \ref{section:mainproof}, we show that at the cost of a constant factor, we can assume the arms remain at least a fixed distance away from each other. In the final part of Section \ref{section:mainproof}, we prove the main result.

\subsection{Acknowledgement} We thank Jack Hanson for useful comments on a draft of this note.

\section{A Shifting Transformation}\label{section:shift-transformation}
We first show the following lemma which is the main ingredient of the proof of the main result. We introduce a transformation that shifts a configuration in a region by $(1/2,1/2)$ to convert arms between the primal and the dual lattices.

Let us begin defining what is meant by a \emph{region} of an annulus bounded by two curves. For this, let $n<N$ and let $\gamma_1$ and $\gamma_2$ be two (primal or dual) paths connecting $\partial B_n$ and $\partial B_N$. In particular, if $\gamma_1$ or $\gamma_2$ is a dual path, we use ``connecting'' in the sense as in the arm events above, which is distinct from the usual topological sense. Consider $\gamma_1$, $\gamma_2$, as well as $\partial B_n$, $\partial B_N$ as curves in $\R^2$. We define $\phi(\gamma_1, \gamma_2)\subset B(n,N)$ to be the Jordan curve\footnote{A Jordan curve is a simple closed curve in $\mathbb{R}^2$.} obtained by concatenating $\gamma_1$, a portion of $\partial B_N$, $\gamma_2$, and a portion of $\partial B_n$ in counterclockwise orientation. If $\gamma_1$ (similarly, $\gamma_2$) lies on the dual lattice and does not (topologically) connect $\partial B_n$ and $\partial B_N$, we add the shortest line segment connecting the endpoints of $\gamma_1$ to $\partial B_n$ or $\partial B_N$ to ensure we obtain a closed curve. 

\begin{definition}
    A \emph{region} is a connected set of edges. The \emph{region in $B(n,N)$ with boundary $\phi$} is the set $S$ of edges whose interiors lie in the interior of $\phi$, together with all edges of $(\phi \cap \partial B_n)\cup (\phi \cap \partial B_N)$.
\end{definition}

\begin{definition}
    For a region $S\subset \mathcal{E}$, we define the event $A_{k,\sigma}(S)$ in a similar manner to the arm event $A_{k,\sigma}(n,N)$, with the additional condition that the $k$ arms consist of edges of $S$ or edges dual to $S$.
    The $k$ arms in $A_{k,\sigma}(S)$ are automatically disjoint from $\partial S\setminus (\partial B(n,N))$. 
\end{definition}

We now introduce a variant of arm events $A_{k,\sigma}$ which includes a separation condition. We fix some integer constant $\ell \ge 5$, the separation between the arms.

\begin{definition}\label{def: separated}
    Let $\tilde{A}_{k, \sigma}(n,N)$ be the event that there are $k$ disjoint arms from distance $n$ to $N$, color-coded by $\sigma$ and any two of the arms are at distance at least $\ell$ between distances $2n$ and $N/2$.
    
    Similarly, for a region $S$ with boundary $\phi$, we define $\tilde{A}_{k,\sigma}(S)$ as the event that there are $k$ disjoint arms connecting $\partial B_n$ and $\partial B_N$ in $S$, color-coded by $\sigma$, and the $k$ arms are at distance at least $\ell$ from each other and $\phi$.
\end{definition}

\begin{lemma} \label{lemma:regionshift}
    Let $k\geq 2$ and $\sigma$ be some color sequence of length $k$. For $n_0(k)<n<N$, let $\phi = \phi(\gamma_1, \gamma_2)$ be a Jordan curve given by two disjoint arms $\gamma_1$, $\gamma_2$, and $S\subset \mathcal{E}$ be the region with boundary $\phi$, excluding any edges in $\gamma_1\cap \gamma_2$. Then,
    \begin{equation} 
        \P(\tilde{A}_{k,\sigma}(S)) \le \P\big(A_{k,\sigma^*}(S\cap B(2n, N/2))\big).
    \end{equation}
\end{lemma}

\begin{proof}
    Our goal is to define an invertible (and thus measure preserving) transformation $T$ on configurations of the edges in $S$ such that
    \begin{equation}\label{eqn: T-sub}
        T(\tilde{A}_{k,\sigma}(S))\subset A_{k,\sigma^*}(S\cap B(2n, N/2)).
    \end{equation}
    Since all configurations are uniformly distributed, we have
    \begin{equation*}
        \P(\tilde{A}_{k,\sigma}(S))
        =\sum_{\omega\in \tilde{A}_{k,\sigma}(S)}\mathbb{P}(\omega) =\sum_{\omega \in \tilde{A}_{k,\sigma}(S)} \mathbb{P}(T(\omega))
    \end{equation*}
    Using the bijectivity of $T$ and \eqref{eqn: T-sub} respectively, the above is bounded by
    \begin{align*}
        \P(\tilde{A}_{k,\sigma}(S))
        &\le \sum_{\omega'\in T(\tilde{A}_{k,\sigma}(S))}\mathbb{P}(\omega')\\
        &\le \mathbb{P}(A_{k,\sigma^*}(S\cap B(2n, N/2))).
    \end{align*}
    This completes the proof given $T$ as desired.
    
    To define $T$, we first choose some deterministic ordering of all edges in $B_N$. This induces an ordering of the edges in $S$, which we enumerate as $e_1,\dots, e_m$.
    
    Given an initial configuration $(\omega(e_1), \dots, \omega(e_m))\in \{0,1\}^S$, we determine an image configuration \[(\omega'(e_1), \dots, \omega'(e_m))\] as an intermediate step to defining $T$ by the following correspondence:
    
    \begin{itemize}
        \item If $(e_i)^* - (1/2,1/2) \in S$, then we let
        \begin{equation*}
            \omega'^*((e_i)^*) =\omega'(e_i) := \omega\left((e_i)^* - \left(\frac12,\frac12\right)\right).
        \end{equation*}
        (Note that $(e_i)^* - (1/2,1/2)$ is an edge on the primal lattice $\mathcal{E}$.) In this case, we say that $e_i$ inherited its status in $\omega'$ from (the status of) $(e_i)^* - (1/2,1/2)$ (in $\omega$).
        
        \item If $(e_i)^* - (1/2,1/2)\not\in S$, then the status of $e_i$ in $\omega'$ remains the same as in $\omega$:
        \[\omega'(e_i) := \omega(e_i)\]
        In this case, we say $e_i$ inherited its status in $\omega'$ from itself.
    \end{itemize}
    See Figure \ref{subfig: inheritance} for an illustration.

    \begin{figure}
        \centering
        \begin{subfigure}[htbp]{.5\textwidth}
            \centering\includegraphics[width=.7\linewidth]{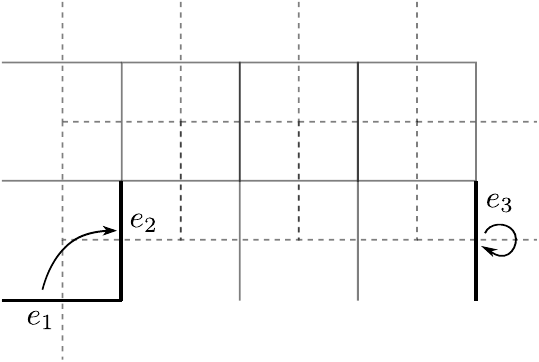}
            \caption{$e_2$ inherits its status in $\omega'$ from $e_1$ (in $\omega$).$\qquad$ \linebreak
                $e_3$ inherits its status in $\omega'$ from itself.}
            \label{subfig: inheritance}
        \end{subfigure}%
        \begin{subfigure}[htbp]{.5\textwidth}
            \centering
            \includegraphics[width=.7\linewidth]{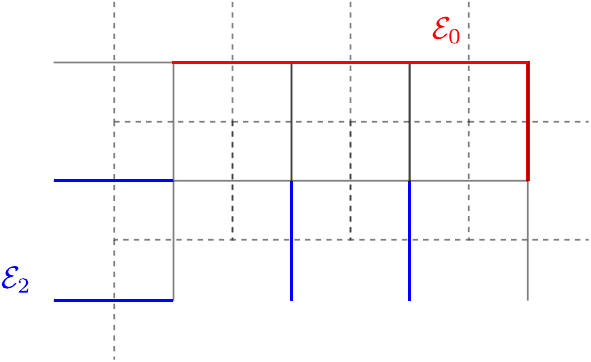}
            \caption{Under transformation $T$, the red edges belong to $\mathcal{E}_0$; the blue edges belong to $\mathcal{E}_2$; and the black edges belong to $\mathcal{E}_1$.}
            \label{subfig: classification}
        \end{subfigure}
        \caption{All solid lines are edges in $S$ and all dotted lines are dual edges of the edges in $S$.}
        \label{fig: inheritance-classification}
    \end{figure}
    
    We classify the edges of $S$ into three sets:
    \begin{enumerate}
        \item An edge $e$ is in $\mathcal{E}_0(T)$ if no edge inherits its status in $\omega'$ from $e$ in $\omega$.
        
        \item An edge $e$ is in $\mathcal{E}_1(T)$ if exactly one edge (including possibly $e$ itself) inherits its status in $\omega'$ from $e$ in $\omega$.
        
        \item An edge $e$ is in $\mathcal{E}_2(T)$ if two edges (including possibly $e$ itself) inherit their status in $\omega'$ from $e$ in $\omega$.
    \end{enumerate}
    Notice that the sets $\mathcal{E}_i$, $i=0,1,2$ do not depend on $\omega$: $e\in \mathcal{E}_0$ if $e^*-(1/2,1/2)\in S$ but $e^*+(1/2,1/2)\notin S$, $e\in \mathcal{E}_2$ if $e^*+(1/2,1/2)\in S$ but $e^*-(1/2,1/2)\notin S$, and an edge $e$ is in $\mathcal{E}_1$ if either  $e^*\pm (1/2,1/2)\in S$ or
    $e^*\pm (1/2,1/2)\notin S$.
    
    By counting the number of edges inheriting their status from each of the sets $\mathcal{E}_i$, $i=0,1,2$, we have: 
    \begin{equation*}
        |S|=|\mathcal{E}_0|+|\mathcal{E}_1|+|\mathcal{E}_2| = 0\cdot |\mathcal{E}_0| + |\mathcal{E}_1| + 2|\mathcal{E}_2|
    \end{equation*}		
    Therefore, $|\mathcal{E}_0|=|\mathcal{E}_2|$. See Figure \ref{subfig: classification} for an illustration.
    
    We now assign new statuses to the edges in $\mathcal{E}_2$. Enumerate the edges in $\mathcal{E}_0$ and $\mathcal{E}_2$ according to the deterministic order fixed in the beginning so that 
    \begin{align*}
        \mathcal{E}_0 &= (e_{\alpha_1}, \dots, e_{\alpha_K}), \\
        \mathcal{E}_2 &= (e_{\beta_1}, \dots, e_{\beta_K})
    \end{align*}
    where $K = |\mathcal{E}_0| = |\mathcal{E}_2|$. Set 
    \begin{equation} \label{eqn:2ndstep}
        T(\omega)(e_{\beta_i}) :=\omega(e_{\alpha_i}),\quad i=1,\dots, K.
    \end{equation}
    For the remaining edges $e\notin \mathcal{E}_2$, we let
    \[T(\omega)(e):=\omega'(e).\]
    

    This definition guarantees that $T$ is invertible. It is easily checked that its inverse is the following map $T'$: given an initial configuration $(\omega'(e_1), \dots, \omega'(e_m))$, first assign the status of each edge $e_i$ as follows.
    \begin{equation}\label{eqn: pre-T'}
        \omega(e_i)=\begin{cases} \omega'((e_i)^* + (\frac12,\frac12)) & \text{ if $(e_i)^* + (\frac12,\frac12) \in S$,}\\
            \omega(e_i) & \text{otherwise}.
        \end{cases}
    \end{equation}
    Next, define $\mathcal{E}_0' = \mathcal{E}_0(T')$ to be the set of edges $e$ such that no edge inherits its status in $\omega$ from $e$, and $\mathcal{E}_2' = \mathcal{E}_2(T')$ to be the set of edges such that two edges inherit their status in $\omega$ from $e$. Here, ``inheritance'' is defined as it was used in the definition of the transformation $T$. Then,
    \begin{align*} 
        \mathcal{E}_2'&= \mathcal{E}_0(T),\\
        \mathcal{E}_0'&= \mathcal{E}_2(T).
    \end{align*}
    For $e\notin \mathcal{E}'_2$, we define
    \begin{equation*}
        T'(\omega')(e)=\omega(e),
    \end{equation*}
    with $\omega$ as in \eqref{eqn: pre-T'}. For $e\in \mathcal{E}_2'=\mathcal{E}_0(T)=(e_{\alpha_1},\ldots,e_{\alpha_K})$, we let
    \begin{equation}\label{eqn: T'}
        T'(\omega'(e_{\alpha_i}))=\omega'(e_{\beta_i}).
    \end{equation}

    Now, we show \eqref{eqn: T-sub}. For any arm $\gamma$ in a configuration $\omega \in \tilde A_{k,\sigma} (S)$, the translated path $\gamma^*$ is contained in $S$ at least from $\partial B_{2n}$ to $\partial B_{N/2}$, since the original arms are $\ell$-separated in the sub-annulus $B(2n,N/2)$. 
    
    Finally, $\gamma^*\subset B(n+1,N-1)$ in the configuration $T(\omega)$ receives the same color as $\gamma$ in the configuration $\omega$ since all but two edges (the two end edges with endpoints lying on $\partial B_n$ or $\partial B_N$) in $\gamma^*$ lie in the set $\mathcal{E}_1$, and thus receive their status from the edges in $\omega$.
\end{proof}


\section{Proof of Proposition \ref{lemma:mainresult}} \label{section:mainproof}

To apply Lemma \ref{lemma:regionshift} in the proof of the main result, we need the arms in the shifted region to be at least distance $\ell$ apart. This is well known to happen with high probability, since if any two arms come close, they form a ``bottleneck'' that generates a $6$-arm event. Moreover, each additional arm that comes close generates two more arms. For completeness, we provide the details in the next lemma.

For $\delta>0$, let $m = m(\delta)$ be the least integer such that $(\pi_1(n))^{2m+2} > n^{-2-\delta}$. For simplicity of notation, we also define $L = 2^{\lfloor \log(d(e)/2)\rfloor}$ \footnote{$\log$ used in this paper denotes logarithm with base $2$.}.
\begin{lemma}\label{lemma:bottleneck}
    Let $d(e)= \mathrm{dist}(e,\partial B_n)$, $k\ge 4$ and $8\ell n_0(k)\leq 8n\leq N$. Suppose $A_{k,\sigma}\setminus \tilde{A}_{k,\sigma}$ occurs, where $\sigma = (\sigma_1, \dots, \sigma_k)$. There exist two arms $\gamma_{-1}$ and $\gamma_{1}$ such that the shortest distance between edges of $\gamma_1$ and $\gamma_{-1}$ in $B(2n,N/2)$ is less than $\ell$ and $\gamma_{-1}$ is of color $\sigma_{i}$, $\gamma_1$ is of color $\sigma_{i+1}$ for some $1\le i\le k$. We relabel the colors $\sigma_{i}$ and $\sigma_{i+1}$ as $\beta_{-1}$ and $\beta_1$, respectively.
    
    Then there exists an edge $e\in B(2n,N/2)$ and integer scales $\lceil\log \ell \rceil= \ell_0 \le \ell_1 \le \cdots \le \ell_{m-1} \le \log(d(e)/2)$, such that the following event $E(e)$ occurs:
    \begin{enumerate}
        \item For $0\le i \le m-2$, if $\ell_{i+1} \ge \ell_{i} + 2$, there are $6+2i$ disjoint arms from $\partial B(e, 2^{1+\ell_{i}})$ to $\partial B(e, 2^{\ell_{i+1}})$. There exists an integer $0\le j\le i$ such that the $6+2i$ arms appear with the color sequence $\Sigma_i$ given by:
        \begin{equation*}
        \begin{split}
            \beta_{-1},\ \beta_{-2}, \dots, \beta_{-i+j-1},\ &\overline{\beta_{-i+j-1}},\ \beta_{-i+j-1}, \dots \\ 
            &\beta_{-2},\ \beta_{-1},\ \beta_1,\ \beta_2, \dots,\ \beta_{j+1},\ \overline{\beta_{j+1}},\ \beta_{j+1}, \dots,\ \beta_2,\ \beta_1. 
        \end{split}
        \end{equation*}

        \item If $\ell_{m-1} + 1 < \log (d(e)/2)$, there are $2m+2$ disjoint arms from $\partial B(e, 2^{1+\ell_{m-1}})$ to $\partial B(e, L)$. There exists an integer $0 \le j\le m-1$ such that the $2m+2$ arms are of the color sequence $\Sigma_{m-1}$:
        \begin{equation*}
            \beta_{-1},\ \beta_{-2}, \dots,\ \beta_{-m+j},\ \beta_{-m+j}, \dots,\ \beta_{-2},\ \beta_{-1},\ \beta_1,\ \beta_2, \dots,\ \beta_{j+1},\ \beta_{j+1}, \dots,\ \beta_2,\ \beta_1.
        \end{equation*}
    \end{enumerate}
\end{lemma}

\begin{proof}
    We first discuss how to determine $\gamma_{-1}$ and $\gamma_1$. On the event that $A_{k,\sigma}$ occurs but $\tilde{A}_{k,\sigma}$ does not, there are $k$ disjoint arms with the color sequence $\sigma$ and there exists a pair of two arms such that the shortest distance between them in $B(2n,N/2)$ is less than $\ell$ and they are of colors $\sigma_i$ and $\sigma_{i+1}$ for some $i$ between $1$ and $k$. If several choices are available, choose the first pair of arms in some fixed deterministic order. We denote the arm of color $\sigma_i$ by $\tilde \gamma_{-1}$ and the arm of color $\sigma_{i+1}$ by $\tilde \gamma_{1}$. Since the two arms come close, the set $\mathcal E \subset B(2n,N/2)$ of edges $e$ such that $B(e,\ell)$ intersects both $\tilde \gamma_{-1}$ and $\tilde\gamma_{1}$ is non-empty. We choose $e\in \mathcal E$ such that $\mathrm{dist}(e,\tilde\gamma_{-1})+ \mathrm{dist}(e,\tilde\gamma_{1})$ is minimized. 
    We then choose $\gamma_{-1}$ and $\gamma_{1}$, two arms of colors $\sigma_i$ and $\sigma_{i+1}$ such that both arms intersect $B(e,\ell)$ and the area enclosed by $\gamma_{-1}, \gamma_{1}, \partial B_n$ and $\partial B_N$ containing $e$ is maximized. We assume that this region lies clockwise of $\gamma_{-1}$ inside $B(n,N)$ and counterclockwise of $\gamma_1$ (otherwise just exchange their labels for the rest of the proof).

    From now on we relabel the colors $\sigma_i$ and $\sigma_{i+1}$ as $\beta_{-1}$ and $\beta_1$. Then, there are four disjoint arms (two of color $\beta_{-1}$ and two of color $\beta_1$) from $B(e, \ell)$ to $B(e, d(e)/2)$ following $\gamma_{-1}$ and $\gamma_1$. See Figure \ref{fig:six-arm}.
    
    \begin{figure}[htbp]
        \centering
        \includegraphics[width=0.45\textwidth]{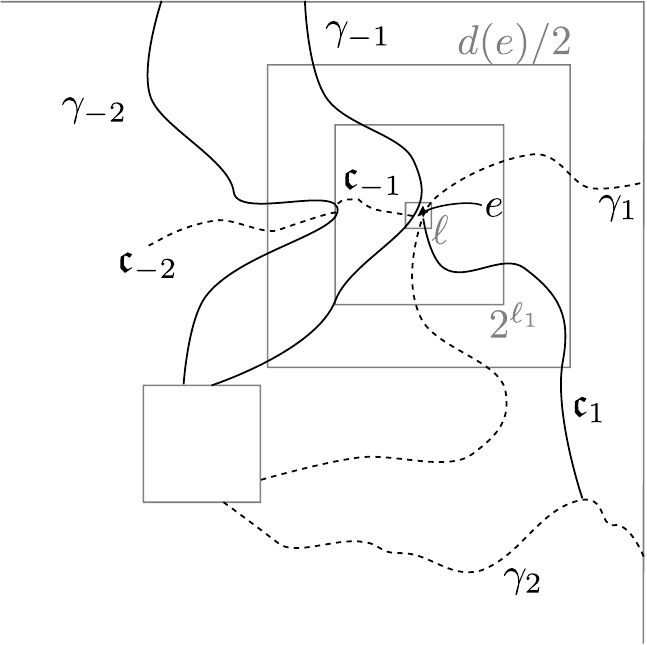}
        \caption{A representation of the arm construction in the proof of Lemma \ref{lemma:bottleneck}. The outer grey boundary represents a portion of $\partial B_N$; the small grey box represents $\partial B_n$.}
        \label{fig:six-arm}
    \end{figure}
    
    Let $\gamma_{-2}$ be the open or dual-closed arm disjoint from $\gamma_{-1}$ such that the area enclosed by $\gamma_{-1}, \gamma_{-2}, \partial B_n$, and $\partial B_N$ on $\gamma_{-1}$'s counterclockwise side is minimized, that is, $\gamma_{-2}$ is the ``closest'' disjoint arm on $\gamma_{-1}$'s counterclockwise side. We denote the color of $\gamma_{-2}$ by $\beta_{-2}$. Similarly, let $\gamma_2$ denote the open or dual-closed arm disjoint from $\gamma_1$ such that the area enclosed by $\gamma_1, \gamma_2, \partial B_n$, and $\partial B_N$ on $\gamma_1$'s clockwise side is minimized. We denote the color of $\gamma_2$ by $\beta_2$. By duality and the minimality of $\gamma_{-2}$, there is a $\overline{\beta_{-1}}$ arm, $\mathfrak{c}_{-1}$, from an edge of $\gamma_{-1}$ in $B(e,\ell)$ to $\gamma_{-2}$. Similarly, there is a $\overline{\beta_1}$ arm, $\mathfrak{c}_1$, from an edge of $\gamma_1$ in $B(e,\ell)$ to $\gamma_2$. See Figure \ref{fig:six-arm}.

    Let $\ell' \ge \log \ell$ be the largest integer such that neither $\gamma_{-2}$ nor $\gamma_2$ intersects $B(e, 2^{\ell'})$. If $2^{\ell'} \ge L$, we let $\ell_1 = \lfloor \log(d(e)/2)\rfloor$. Then there are six arms in the annulus $B(e,\ell, d(e)/2)$: a portion of $\gamma_1$, a portion of $\mathfrak{c}_1$, another portion of $\gamma_1$, a portion of $\gamma_{-1}$, a portion of $\mathfrak{c}_{-1}$, and another portion of $\gamma_{-1}$. Otherwise if $2^{\ell'} < L$, we let $\ell_1 = \ell'$ and we have the same arms as in the previous case crossing $B(e, \ell, 2^{\ell_1})$ instead. In this case, without loss of generality, we assume $\gamma_{-2}$ intersects $B(e, 2^{\ell_1 + 1})$. We thus associate to the scale $\ell_1$ the collection
    \[\mathcal{C}_1=\{\gamma_{-2}, \gamma_{-1},\gamma_1\}\]
    of three arms intersecting $B(e, 2^{\ell_1 +1})$.

    For any $1\le i\le m-2$, we inductively define the scale $\ell_{i+1}$. If $\ell_{i} = \lfloor\log (d(e)/2)\rfloor$, we let all subsequence $\ell_{r} = \ell_{i}$ for $r\ge i+1$. Otherwise, $\ell_{i} < \lfloor \log (d(e)/2)\rfloor$ and there is a collection \[\mathcal{C}_{i}=\{\gamma_{-i+j-1}, \dots, \gamma_{-1}, \gamma_1, \dots, \gamma_{j+1}\}\]
    of $i+2$ arms associated with the scale $\ell_{i}$, with $j$ ($j\le i$) arms on the clockwise side of $\gamma_1$. These arms all intersect $B(e, 2\cdot 2^{\ell_{i}})$. Their colors are labeled $\beta_{-i+j-1}, \dots, \beta_{-1}, \beta_1, \dots, \beta_{j+1}$. Let $\gamma_{-i+j-2}$ be the open or dual-closed arm disjoint from $\gamma_{-i+j-1}$ such that the area of the region enclosed by $\gamma_{-i+j-1}, \gamma_{-i+j-2}, \partial B_n$, and $\partial B_N$ on $\gamma_{-i+j-1}$'s counterclockwise side is minimized. Its color is labeled $\beta_{-i+j-2}$. Similarly, we define $\gamma_{j+2}$ the ``closest'' arm to the clockwise side of $\gamma_{j+1}$ and label its color $\beta_{j+2}$.

    By duality and the minimality of $\gamma_{-i+j-2}$, there is a $\overline{\beta_{-i+j-1}}$ arm, $\mathfrak{c}_{-i+j-1}$, from an edge of $\gamma_{-i+j-1}$ in $B(e, 2^{1+\ell_{i-1}})$ to $\gamma_{-i+j-2}$. Similarly, there is a $\overline{\beta_{j+1}}$ arm, $\mathfrak{c}_{j+1}$, from an edge of $\gamma_{j+1}$ in $B(e, 2^{1+\ell_{i}})$ to $\gamma_{j+2}$.

    Let $\ell'' \ge \ell_{i}$ be the largest integer such that neither $\gamma_{-i+j-2}$ nor $\gamma_{j+2}$ intersects $B(e,2^{\ell''})$. If $2^{\ell''} \ge L$, we let $\ell_{i+1} = \lfloor\log (d(e)/2)\rfloor$. Otherwise, we let $\ell_{i+1} = \ell''$. As in the $i=1$ case, we form $\mathcal{C}_{i+1}$ from $\mathcal{C}_{i}$ by adding the arm, either $\gamma_{-i+j-2}$ or $\gamma_{j+2}$, that intersects $B(e,2^{1+\ell_{i+1}})$ for this scale. Suppose, without loss of generality, that the newly recorded arm is $\gamma_{j+2}$:
    \[\mathcal{C}_{i+1}:=\mathcal{C}_{i}\cup \{\gamma_{j+2}\}.\]
    The $i+3$ arms in $\mathcal{C}_{i+1}$ cross the box $B(e, 2^{1+\ell_{i+1}})$. We associate the collection $\mathcal{C}_{i+1}$ to the scale $\ell_{i+1}$. Each of the arms in $\mathcal{C}_{i}$ crosses the annulus $B(e, 2^{1+\ell_{i}}, 2^{\ell_{i+1}})$ twice. In addition, both $\mathfrak{c}_{-i+j-1}$ and $\mathfrak{c}_j$ cross $B(e, 2^{1+\ell_{i}}, 2^{\ell_{i+1}})$. Thus, the annulus is crossed by $6+2i$ arms.
    
    For the $m$-th annulus, if $\ell_{m-1} + 1 < \lfloor \log(d(e)/2)\rfloor$, then each of the $m+1$ arms in $\mathcal C_{m-1}$ crosses $B(e, 2^{1+\ell_{m-1}}, L)$ twice, proving item 2.
    
\begin{figure}[htbp]
    \centering
    \includegraphics[width=0.6\textwidth]{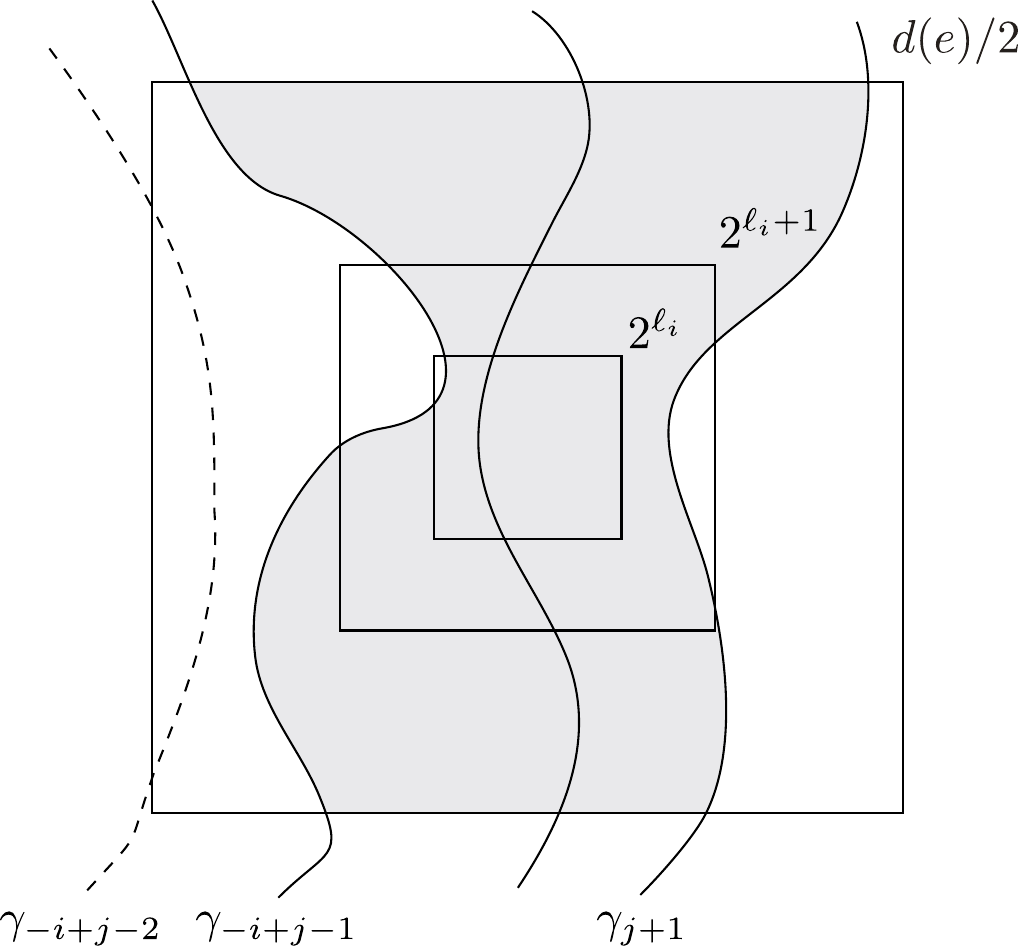}
    \caption{The geometric configuration in the proof of disjointness of Lemma \ref{lemma:bottleneck}. The shaded region represents $R_2$.}
    \label{fig:disjointness}
\end{figure}

    Finally, we address the disjointness of arms in each annulus. For the first annulus $B(e, \ell, 2^{\ell_1})$, $\gamma_{-1}$, $\gamma_{1}$, $\mathfrak{c}_{-1}$, and $\mathfrak{c}_1$ are disjoint between $\partial B(e, \ell)$ and $\partial B(e, 2^{\ell_1})$ by definition. Inductively, for any non-empty annulus $B(e, 2^{1+\ell_{i}}, 2^{\ell_{i+1}})$, $0\le i \le m-2$, we have associated a collection $\mathcal{C}_{i}$ of $i+2$ arms to the scale $\ell_{i}$. Suppose all of them are disjoint. The ``outermost'' arms, labeled $\gamma_{-i+j-1}$ and $\gamma_{j+1}$, partition the box $B(e, 2^{\ell_{i+1}})$ into three regions: $R_1$, the region on the counterclockwise side of $\gamma_{-i+j-1}$; $R_2$, the region between $\gamma_{-i+j-1}$ and $\gamma_j$ through which the other $i$ arms; and $R_3$, the region on the clockwise side of $\gamma_{j+1}$. The collection $\mathcal{C}_{i+1}$ is formed by adding either $\gamma_{-i+j-2}$ or $\gamma_{j+2}$ to $\mathcal{C}_{i}$. Neither passes through $R_2$ and both are necessarily disjoint from $\gamma_{-i+j-1}$ and $\gamma_{j+1}$. Thus they must be disjoint from all $i+2$ arms in $\mathcal C_i$. Thus, all $i+3$ arms in $\mathcal C_{i+1}$ are disjoint. See Figure \ref{fig:disjointness} for an illustration.
\end{proof}

In the previous proof, $B(e, \ell, L)$ was divided into $m$ annuli containing six arms in $B(e, \ell, 2^{\ell_1})$ and two additional arms for every successive scale larger than $\ell_1$. Note that the $\ell_i$, $i\le r \le m$ might coincide for some $r<m$ if we reach distance $L$ before the $m$-th step of induction. In this case, the corresponding annuli $B(e,2^{1+\ell_i}, 2^{ \ell_{i+1}})$ are empty.

\begin{corollary} \label{cor:non-consecutive}
    For $0\le i \le m-2$, each $(6+2i)$-arm event is of a color sequence that contains at least two non-consecutive occurrences of open colors and resp. dual-closed colors.
\end{corollary}

\begin{proof}
    We recall that for $0\le i\le m-2$, the $(6+2i)$-arm color sequence is of the form
    \begin{equation*}
        \begin{split}
            \beta_{-1},\ \beta_{-2}, \dots, \beta_{-i+j-1},\ &\overline{\beta_{-i+j-1}},\ \beta_{-i+j-1}, \dots \\ 
            &\beta_{-2},\ \beta_{-1},\ \beta_1,\ \beta_2, \dots,\ \beta_{j+1},\ \overline{\beta_{j+1}},\ \beta_{j+1}, \dots,\ \beta_2,\ \beta_1. 
        \end{split}
        \end{equation*}

    If $\beta_{j+1} = \beta_{-i+j-1}$, then the two occurrences of $\beta_{j+1}$ are non-consecutive and $\overline{\beta_{j+1}}$ and $\overline{\beta_{-i+j-1}}$ are two non-consecutive occurrences of the $\overline{\beta_{j+1}}$ color.

    If $\beta_{j+1} = \overline{\beta_{-i+j-1}}$, then the two occurrences of $\beta_{j+1}$ are non-consecutive and $\overline{\beta_{j+1}}$ and $\beta_{-i+j-1}$ are two non-consecutive occurrences of the $\overline{\beta_{j+1}}$ color.
\end{proof}

A landing sequence $I = \{I_i\}_{1\leq i\leq k}$ on $\partial B_n$ is a sequence of disjoint sub-intervals of $\partial B_n$ in clockwise order where each $|I_i| \ge \delta n$ for some $\delta>0$. A landing sequence $I'$ on $\partial B_N$ is defined analogously. We use a result in \cite{Nolin08}, where Nolin proved that the probability of a $k$-arm event is comparable to the probability of the same event with extra landing conditions. 
\begin{equation} \label{eqn:landing}
    \P(A_{k,\sigma}(n,N))\asymp \P(A_{k,\sigma}^{I/ I'}(n,N)),
\end{equation}
see \cite[Theorem 11]{Nolin08}. The theorem is originally stated with additional ``well-separatedness'' information which is omitted here.

Even though the setting of Nolin's paper is site percolation on the triangular lattice, the exact proof applies to the square lattice as the main techniques used in the proof are Russo-Seymour-Welsh estimates and generalized FKG inequality, both of which apply to the square lattice. Using \eqref{eqn:landing}, we may work with the events with prescribed landing zones for all arms.

The next lemma shows that enforcing a separation condition in the arm events does not essentially change the order of the probabilities.
\begin{lemma}\label{lemma:separatedness}
    For Bernoulli percolation on the square lattice, any integer $k\geq 2$, some color sequence $\sigma$, two landing sequences $I, I'$ on $\partial B_n$ and $\partial B_N$, and any $\epsilon>0$. There exists $n$ sufficiently large and $N \ge 8n$ such that
    \begin{equation}\label{eq:wellseparatedness}
        \P(\tilde{A}_{k,\sigma}^{I,I'}(n,N)) \geq (1-\epsilon) \P(A_{k,\sigma}^{I/I'}(n,N)).
    \end{equation}
\end{lemma}

\begin{proof}
    Equivalently, we show the following:
    \begin{equation*}
        \P(A^{\times, I, I'}_{k,\sigma}(n,N)) \leq \epsilon \P(A_{k,\sigma}^{I,I'}(n,N)),
    \end{equation*}
    where $A^{\times, I,I'}_{k,\sigma}$ denotes the event $A_{k,\sigma}^{I,I'} \setminus \tilde{A}^{I,I'}_{k,\sigma}$: that is, there are two arms that come closer than $\ell$ in $B(2n,N/2)$ in a $k$-arm event. 
    
    Let the event $E(e)$ be as in Lemma \ref{lemma:bottleneck}, then
    \begin{equation*}
        \P(A^{\times, I,I'}_{k,\sigma}(n,N))\leq \sum_{e\in B(2n,N/2)} \P(A^{\times,I,I'}_{k,\sigma}(n,N), E(e))
    \end{equation*}
    
    Using independence, Lemma \ref{lemma:bottleneck}, and carving out a sub-annulus to distance $d(e)/2$ around $e$, the terms on the RHS can be bounded as
    \begin{align}
         \P&(A^{\times,I,I'}_{k,\sigma}(n,N), E(e)) \nonumber\\
            &\quad \leq \P(A^{I}_{k,\sigma}(n, d(e)/2)) \P(A^{I'}_{k,\sigma}(2d(e), N)) \nonumber\\
            &\quad \times \sum_{\ell_0\le \ell_1\le \cdots\le \ell_{m-1}}^{\lfloor\log(d(e)/2)\rfloor} \prod_{i=0}^{m-2} \P(A_{6+2i, \Sigma_i}(e, 2^{\ell_i+1}, 2^{\ell_{i+1}})) \P(A_{2m+4, \Sigma_{m-1}}(e, 2^{\ell_{m-1}+1}, L). \label{eq:btn-prod}
    \end{align}
    
    By the generalized FKG inequality and quasi-multiplicativity, see \cite[Proposition 12]{Nolin08}, we have
    \begin{equation*}
        \P(A^{I}_{k,\sigma}(n, d(e)/2)) \P(A^{I'}_{k,\sigma}(2d(e), N)) \leq C\P(A^{I,I'}_{k,\sigma}(n,N))
    \end{equation*}
    where the constant $C$ depends only on $k$ and $\sigma$.
    
    On the other hand, we know that the alternating five-arm probability has a universal exponent $2$, see \cite[Lemma 5]{KSZ98}:
    \begin{equation*}
        \P(A_{5, OC^*OOC^*}(n,N)) \le C \left(\frac{n}{N}\right)^{-2}
    \end{equation*}
    for some constant $C>0$.
    By Corollary \ref{cor:non-consecutive}, for any $0\le i \le m-2$, the color sequence $\Sigma_i$ can be split into an alternating five-color sequence ($(O,C^*,O,O,C^*)$ or its flipped sequence) and a color sequence of length $2i+1$, while maintaining the relative order within each sub-color sequence. We can now apply Reimer's inequality, see \cite{Reimer00}, and obtain
    \begin{equation*}
        \begin{split}
            \P(A_{6+2i, \Sigma_i}(e, 2^{\ell_i+1},2^{\ell_{i+1}}))
            &\le \pi_5(2^{\ell_i+1},2^{\ell_{i+1}}) \P(A_{2i+1}(e, 2^{\ell_i+1},2^{\ell_{i+1}})) \\
            &\le \pi_5(2^{\ell_i+1},2^{\ell_{i+1}}) (\pi_1(2^{\ell_i+1},2^{\ell_{i+1}}))^{2i+1}.
        \end{split}  
    \end{equation*}
    Here we denote $\pi_5(n,N)=\mathbb{P}(A_{5,OC^*OOC^*}(n,N))$ and $\pi_1(n,N)=\mathbb{P}(A_1(n,N))$. For $\pi_1$, we have the estimate
    \begin{equation*}
        \pi_1(n,N) \le \left(\frac{n}{N}\right)^\delta, \quad \text{for some $\delta>0$.}
    \end{equation*}

    Similarly, for the final annulus, we have
    \begin{equation*}
        \P(A_{2m+4, \Sigma_{m-1}}(e, 2^{\ell_{m-1}+1}, L)) \le \left(\frac{2^{\ell_{m-1}}}{L}\right)^{2+(2m-1)\delta}.
    \end{equation*}

    Plugging the above estimates in the product in line \eqref{eq:btn-prod}, we have
    \begin{equation*}
        \begin{split}
            \text{product}
            &\le \prod_{i=0}^{m-2}(2^{\ell_i - \ell_{i+1}})^{2+(2i+1)\delta} \left(\frac{2^{\ell_{m-1}}}{L}\right)^{2+(2m-1)\delta} \\
            &= \left(\frac{\ell}{L}\right)^{2+\delta} \prod_{i=0}^{m-2} (2^{2\delta})^{\ell_i} (2^{2\delta})^{-(m-2)\ell_{m-1}} \left(\frac{2^{\ell_{m-1}}}{L}\right)^{2(m-1)\delta}.
        \end{split}
    \end{equation*}
    Since we have the estimate $\sum_{i\le N} x^i \le c x^N$ for any $x= 2^y$, $y>0$, and some constant $c = c(y)$, we inductively bound 
    \begin{equation*}
        \sum_{\ell_0 \le \ell_i \le \ell_{i+1}} (2^{2i\delta})^{\ell_i} \le c (2^{2i\delta})^{\ell_{i+1}}.
    \end{equation*}
    Thus, summing over all possible values of $\ell_1, \dots, \ell_{m-1}$, we have
    \begin{align*}
        \eqref{eq:btn-prod}
        &\le \left(\frac{\ell}{L}\right)^{2+\delta} c^{m-2} \sum_{ \ell_{m-1} \ge \log(\ell)}^{\lfloor\log(d(e)/2)\rfloor} (2^{2\delta})^{(m-2)\ell_{m-1}} (2^{2\delta})^{-(m-2)\ell_{m-1}}\left(\frac{2^{\ell_{m-1}}}{L}\right)^{2(m-1)\delta} \\
        &\le \left(\frac{\ell}{L}\right)^{2+\delta} c^{m-2} (2^{\lfloor \log(d(e)/2)\rfloor})^{-2(m-1)\delta} c(2^{2(m-1)\delta})^{\lfloor\log(d(e)/2)\rfloor} \\
        &\le C \left(\frac{\ell}{L}\right)^{2+\delta}.
    \end{align*}
    Dyadically summing over the location of $e$, we have
    \begin{align*}
        \P(A^{\times,I,I'}_{k,\sigma}(n,N))
        &\leq C \P(A^{I,I'}_{k,\sigma}(n,N)) \sum_{e\in B(2n,N/2)} (d(e)/2)^{-2-\delta} \\
        &\leq C \P(A^{I,I'}_{k,\sigma}(n,N)) \sum_{k=1}^{\lfloor\log (N/4n)\rfloor} \left(n2^{k-1}
        \right)^{-2-\delta} n^2 2^{2k} \\
        &\leq C' n^{-\delta} \P(A^{I,I'}_{k,\sigma}(n,N)),
    \end{align*}
    where the factor $n^22^{2k}$ estimates the number of edges in the $k$-th annulus of the sum. Choosing $n$ such that $C'n^{-\delta} \leq \epsilon$, the desired result follows.
\end{proof}

Finally, we prove the main result.
\begin{proof}[Proof of Proposition \ref{lemma:mainresult}]
    By Lemma \ref{lemma:separatedness} and \eqref{eqn:landing}, we have for any $\epsilon>0$ and some choice of landing sequences $I,I'$ which we will specify later on, there exists a $C>0$ such that 
    \begin{equation*}
        \P(A_{k,\sigma}(n,N)) \le C\P(A^{I,I'}_{k,\sigma}(n,N))\leq \frac{C}{1-\epsilon}\P(\tilde{A}^{I,I'}_{k,\sigma}(n,N)),
    \end{equation*}
    where $\tilde{A}_{k,\sigma}$ is the $\ell$-separated $k$-arm event, see Definition \ref{def: separated}. It suffices to bound the probability on the RHS.
    
    If suffices to consider the case when $\sigma$ and $\sigma'$ differ by one entry. For general polychromatic color sequences $\sigma$ and $\sigma'$, we consider a sequence interpolating between $\sigma$ and $\sigma'$ with at most $k$ steps such that each two consecutive color sequences differ by one entry.
    
    Without loss of generality, we assume that $\sigma$ and $\sigma'$ differ only in the $k$-th entry and moreover we assume that $\sigma_1$ and $\sigma_2$ are open and dual-closed respectively. Fix two consecutive landing zones on $\partial B_N$ for an open and a dual-closed arm, say $I_1$, $I_2$, corresponding to the first two entries. Let $\gamma_1$, $\gamma_2$ be the pair of open and dual-closed arms closest to each other, such that $\gamma_1$ lands on $I_1$ and $\gamma_2$ lands on $I_2$. Finally, we let $\gamma_3, \dots, \gamma_{k-1}$ be such that $\gamma_i$ is the disjoint arm with color $\sigma_i$ landing on $I_i$ such that the enclosed region between $\gamma_{i-1}$ and $\gamma_i$ is minimal. 
    
    We then denote by $U^c$ the region enclosed by $\gamma_1$, $\gamma_{k-1}$, $\partial B_n$, and $\partial B_N$ that excludes $\gamma_k$. By minimality, the event $\{U^c=R\}$ depends only on the status of edges in $R$. In particular, the configuration in the complement region $U = B(n,N)\setminus U^c$ is independent of $U^c$. Moreover, there is an arm $\gamma_k$ with color $\sigma_k$ in $U$ such that $\gamma_k$ is at distance at least $\ell$ from $\gamma_1$ and $\gamma_{k-1}$, two parts of the boundary of $U$. 
    \begin{align*}
        \P(\tilde{A}^{I_1, \dots, I_k}_{k,O, C^*, \sigma_3, \dots, \sigma_{k}}(n,N)) 
        &\leq \sum_{\text{admissible } S}\P( \tilde{A}^{I_k}_{1,\sigma_{k}}(U), U=S) \\
        &=\sum_{\text{admissible } S} \P(U=S) \P(\tilde{A}^{I_k}_{1,\sigma_{k}}(S))
    \end{align*}
    Each of the sums above is over the possible values of $S$ of (random) region $U$.
    
    With fixed $\gamma_1,\dots, \gamma_{k-1}$, we flip the percolation configuration in the region $S$. We have
    \begin{equation}\label{eqn:flip}
        \P(\tilde{A}^{I_k}_{1,\sigma_{k}}(S)) = \P(\tilde{A}^{I_k}_{1,\overline{\sigma}_{k}}(S)).
    \end{equation}
    Recall that $\overline{\sigma}_k$ denotes the flipped color sequence.
    
    Having flipped the configuration in $S$, we use the transformation $T$ defined in the proof of Lemma \ref{lemma:regionshift} to shift the configuration in the region $S$ to the dual lattice. By Lemma \ref{lemma:regionshift}, we have
    \begin{equation*}
        \P(\tilde{A}^{I_k}_{1, \overline{\sigma}_{k}}(S)) \leq \P( A^{I_k}_{1, \overline{\sigma}^*_{k}}(S\cap B(2n,N/2))).
    \end{equation*}
    We note that even though Lemma \ref{lemma:regionshift} is not stated with prescribed landing zones, the shifting transformation does not essentially change where the arm lands. So the inequality holds as is.

    Inserting this inequality into \eqref{eqn:flip} we have
    \begin{align*}
        \P(A^{I_1, \dots, I_k}_{k,O, C^*, \sigma_3, \dots, \sigma_{k}}(n,N))
        &\leq C \sum_{\text{admissible } S} \P(U=S) \P( A^{I_k}_{1, \overline{\sigma}^*_{k}}(S\cap B(2n,N/2))) \\
        &\leq C\P(A^{I_1, \dots, I_k}_{k,O,C^*,\sigma_3,\dots, \sigma_{k-1}, \overline{\sigma}^*_{k}}(2n,N/2))\\
        &\leq C\P(A_{k,O,C^*,\sigma_3,\dots, \sigma_{k-1}, \overline{\sigma}^*_{k}}(n,N)).
    \end{align*}	
    The final inequality is given by a standard RSW-type argument, see \cite[Proposition 16]{Nolin08}.
    
    %

\end{proof}

\bibliography{reference}
\bibliographystyle{amsplain}

\end{document}